\documentclass{article}
\usepackage{amsmath,amssymb}

\def\L{\mathbb L}
\def\R{\mathbb R}
\def\C{\mathbb C}
\def\W{\mathbb W}

\def\ov{\overline}
\def\wt{\widetilde}
\def\wht{\widehat}
\def\H{\mathcal H}

\def\LL{\mathcal L}
\def\F{\mathcal F}
\def\P{\mathcal P}

\def\epsilon{\varepsilon}
\def\phi{\varphi}
\def\le{\leqslant}
\def\ge{\geqslant}

\renewcommand{\Re}{\mathop{\mathrm{Re}}\nolimits}
\renewcommand{\Im}{\mathop{\mathrm{Im}}\nolimits}
\newcommand{\res}{\mathop{\mathrm{res}}\nolimits}
\newcommand{\sgn}{\mathop{\mathrm{sgn}}\nolimits}

 \newcounter{sec}

\begin{document}

\begin{center}
{\Large \bf

Some remarks on stable densities and
operators of
fractional differentiation
}

\medskip

\large\sc Neretin Yuri A.

\end{center}

\begin{flushright}
\small\it Dedicated to A.M.Vershik
\end{flushright}

{\small  Let $D(s)$ be a fractional derivation of order $s$.
For  real $\alpha\ne 0$,
we construct an integral operator  $A(\alpha)$
in an appropriate functional space such that
$A(\alpha) D(s) A(\alpha)^{-1}=D(\alpha s)$ for all $s$.
The kernel of the operator $A(\alpha)$ is expressed in terms
of a function similar to the stable densities.}

\bigskip

{\bf 0.1. Definition of functions $\L_{\alpha,\beta}$.}
This paper contains several simple observations
concerning the special function
\begin{equation}
\L_{\alpha,\beta}(z)=
\sum_{n=0}^\infty
\frac{ (-1)^n \Gamma(\alpha n+\beta)}
       {n!}
        z^n
\end{equation}
where
$0<\alpha<1$, $\Re \beta>0$,
  and $z\in\C$. We can also represent this function in the
form
\begin{align}
\L_{\alpha,\beta}(z)&=
\int_0^\infty x^{\beta-1} \exp(-z t^\alpha-t)\,dt\\
\L_{\alpha,\beta}(z)&=
 \frac 1{2\pi i}
\int_{-i\infty}^{+i\infty}
\Gamma(\beta-\alpha s) \, \Gamma(s)z^{-s}\,ds
\end{align}
The integrals (0.2), (0.3) also make sense for $\alpha>1$.
The definition of functions $\L_{\alpha,\beta}$
is discussed in details below in Section 1.

The function $\L_{\alpha,\beta}$ is one of the simplest
examples of the so-called {\it $H$-functions}
 (or {\it Fox functions}),
see \cite{Mar}.
In a strange way, the function $\L_{\alpha,\beta}$
 has no official name.
Obviously, for rational $\alpha=p/q$ the function
$\L_{\alpha,\beta}$ can be expressed
 in the terms of higher hypergeometric
functions. But for $q>4$ such expressions
do not seem very useful.

\smallskip

{\bf 0.2. Results of the paper.
Integral operators with functions
$\L_{\alpha,\beta}$ in kernels.}

{\bf A)} We consider the space  $\mathcal K$
of functions
holomorphic in the half-plane $\Re z>0$,
smooth up to the line $\Re z=0$,
 and satisfying the following condition

-- for each $k>0$ and $N>0$ there exists $M$ such that
\begin{equation}
|f^{(k)}(z)|\le M (1+|z|)^{-N}
\end{equation}

We define the operators of
{\it fractional differentiation}
 $D_h$ in the space $\mathcal
K$
by
$$
D_h f(z)=
\frac{\Gamma(h+1)}{2\pi }
\int_{-\infty}^{+\infty}
\frac{f(it)\,dt}{( -it+z)^{h+1}}
$$
For a positive integer $n$, we have $D_n=(-1)^n d^n/ dz^n$;
the operator $D_{-n}$ is
the indefinite integration iterated $n$ times.
Also $D_{h+r}=D_h D_r$. See Section 2 below for details.

Next, for $\alpha>0$ we define the kernel
$$
K_\alpha(u,v)=
\int_0^\infty \exp(-ux^\alpha-vx)\, dx=
v^{-1} \L_\alpha(u/v^\alpha),
$$
and the operator in the space $\mathcal K$
given by
\begin{equation}
A_\alpha f(v)=\frac 1 {2\pi}
\int_{-\infty}^{+\infty}
 K_\alpha(-it,v) f(it)\, dt
\end{equation}
The operators $A_\alpha$ form an one-parameter
group (see Subsection 3.5),
\begin{equation}
A_\alpha A_\beta=A_{\alpha\beta};
\end{equation}
We also show that they satisfy the property
\begin{equation}
A_\alpha D_h  A_\alpha^{-1}=D_{\alpha h}
\end{equation}

{\sc Remark.}
Emphasis the following particular cases of (0.7):
\begin{align*}
A_{m/n} \frac {d^n}{dz^n}A_{m/n}^{-1}&=(-1)^{m-n} \frac {d^m}{dz^m};\\
A_k \frac d{dz} A_k^{-1}&=(-1)^{k-1} \frac {d^k}{dz^k}
\end{align*}
for integer $m$, $n$, $k$.

{\sc Remark.} Also
\begin{equation}
A_\alpha z A_\alpha^{-1} f(z)= \frac 1\alpha D_{1-\alpha} (z f(z))
\end{equation}

{\sc Remark.} It seems that (0.7), (0.8) and formula
(0.12) below give a possibility to strange transformations
of partial differential equations and their solutions.

\smallskip

Further, the operators of dilatation
$$R_a g(z)= a^{-1}g(z/a), \qquad a>0$$
satisfy
\begin{equation}
A_\alpha^{-1} R_a A_\alpha= R_{a^\alpha}
\end{equation}

The generator of the one-parameter group $R_a$ is
$(z\, d/dz+1)$. Hence (0.9) can be written in the form
$$
A_\alpha^{-1} \Bigl(z \frac d{dz}+1\Bigr) A_\alpha
=\alpha\Bigl(z \frac d{dz}+1\Bigr)
$$

{\bf B)} In Section 3, we consider the group $G$ of operators in
$\mathcal K$ generated by the operators $A_\alpha$,
the fractional derivations $D_h$, and the dilatations $R_a$.
We observe that $G$ is a 6-dimensional solvable Lie group with
2-dimensional center and kernels of all elements of this group
admit  simple expressions in the terms
 of the functions $\L_{\alpha,\beta}$
 (Theorems 3.1, 3.2).

{\bf C)}
In Section 4, we consider the usual Riemann--Liouville
fractional integrations $J_h$
in the space of functions on the half-line $x\ge 0$,
 see below (4.1).
For $0<\alpha<1$, we  consider
 the Zolotarev  operators \cite{Zol1}--\cite{Zol2}
defined by the formula
\begin{equation}
B_\alpha f(x)
=\frac 1{\pi x} \int_0^\infty
\Im \Bigl\{\L_{\alpha,1}(x^{-\alpha}y e^{i\pi\alpha})\Bigr\}
f(y)\,dy
\end{equation}
We have
\begin{align}
B_\alpha B_\beta&=B_{\alpha\beta}\\
B_\alpha J_h&=J_{h\alpha} B_\alpha
\end{align}
(but we can not represent the identity (0.12) in the form (0.7),
since the operators $B_\alpha$ are not invertible).

These operators can be included to a 7-dimensional
{\it semigroup} of integral
 operators on the half-line, this semigroup has
   a 2-dimensional center
    (Theorem 4.2).

\smallskip

{\bf 0.3. Some  references on functions $\L_{\alpha,\beta}$.}

1) Barnes in 1906 \cite{Bar}  evaluated asymptotics of
several $H$-functions and, in particular, for
$\L_{\alpha,\beta}$. But it seems that he had no
reasons to investigate $\L_{\alpha,\beta}$ in details;
in the sequel years
 this function
(as far as I know) had not attracted
specialists in special functions.

\smallskip

2) The functions
\begin{equation}
\W_{\alpha,\beta}(z)=
\sum_{n=0}^\infty
\frac{z^n}
{n!\,\Gamma(\alpha n+\beta)},
\qquad \alpha>0
\end{equation}
('{\it Wright functions}', '{\it Bessel--Maitland functions}')%
\footnote{Apparently, Wright and Maitland are coinciding persons
 (Edward Maitland Wright). He also coincides with the author
 of the well-known book Hardy, Wright "An introduction to
 the number theory".}
were discussed more,
see \cite{Wri}, \cite{Aga}, and
 references in "Higher transcendental functions"
\cite{HTF} (Section "Mittag-Leffler function").
 The functions $\L_{\alpha,\beta}(z)$,
$\W_{\alpha,\beta}(z)$ quite often appear in formulae
in similar cases (for instance, see below Subsection 1.5).
Another 'relative' of the function $\L_{\alpha,\beta}$
is the Mittag-Leffler function
$\sum z^n/\Gamma(\alpha n+1)$,
 that appear in literature
quite often.

\smallskip

3) The functions $\L_{\alpha,\beta}$
appear (see Feller \cite{Fel0}) if we
solve the Cauchy problem for
the partial pseudo-differential equation
$$
\Bigl[\frac d{dt} - \frac {d^\alpha}{dx^\alpha} \Bigr]
f(x,t)=0,\qquad f(x,0)=\psi(x)
$$
where
$d^\alpha/dx^\alpha$ is some
fractional derivative.
Sometimes it is possible to write
$$
f(x,t)=\int K(t,x,y) f(y)\,dy
$$
where the kernel $K$ can be expressed in the terms of
the function $\L_{\alpha,\beta}$. The
work of Feller generated a wide literature
on diffusions generated by pseudo-differential operators.

\smallskip

4)
Now, we recall the most important situation, where the functions
 $\L_{\alpha,\beta}$ arise
in a natural way.

Consider a sequence $\xi_j$ of independent random variables
and its partial sums sums $S_n=\xi_1+\dots +\xi_n$.
Consider the distribution $\mu_n$ of $S_n$.
Let us center and normalize $\mu_n$ in some way,
$\wt\mu_n(t):=\mu_n(a_n t +b_n)$, where $a_n>0$, $b_n\in \R$
 are some constants.
Which distributions can appear as limits
of sequences $\wt\mu_n$?
In the most common cases, we  obtain a normal (Gauss) distribution.
Nevertheless, there are other possible limits
(\cite{Lev3}, see also \cite{Fel}), they are named by {\it stable
distributions}.
Densities of these distributions admit
a simple expression (0.15) in terms of the functions
$\L_{\alpha,1}$.

A logical possibility
of non-normal distributions in limit theorems
of this kind
 was  observed by Cauchy
in 1853, see \cite{Cau}. He claimed that the
distribution, whose densities are given by
\begin{equation}
\phi_\alpha(x)=\int_0^\infty \exp(-x^\alpha)\cos(tx)\,dt
\end{equation}
can appear in limit theorems for sums
 of independent random variables.
 Firstly, it was necessary to verify positivity
 of the functions $\phi_\alpha$.
They are is really positive for $0<\alpha\le 2$,
but Cauchy could not prove this
except several simple  cases ($\alpha=1,1/2,2$).
In 1922 P.Levy%
\footnote{Some references between Cauchy and Levy
can be found in \cite{Win}. Also,
there was work of Holtsmark (1919) on distribution
of gravitation force
in Universe, see its exposition in \cite{Fel},
\cite{Zol2}}
 attracted attention to the problem
(\cite{Lev1}), and
 in 1923 Polya \cite{Pol}
  proved positivity of (0.14) for $0<\alpha<1$.

 After appearance of Kolmogorov--Levy--Hinchin
 integral representation
for infinitely divisible laws,
a complete description of stable distributions became
 solvable problem, the final result is present in
the books of P.Levy \cite{Lev2}, 1937, and
A.Hinchin (another spelling is 'Khintchine'),
 \cite{Hin}, 1938.
The stable densities can be represented in the form
\begin{equation}
p(x;\alpha,\gamma)=
\frac 1{\pi x} \Im \L_{\alpha,1}
( x^{-\alpha} e^{i(\gamma-\alpha)\pi/2})
\end{equation}
where $0<\alpha<2$, $\gamma\in\R$, and
 $|\gamma|<\min (\alpha, 2-\alpha)$
(in this formula,
we omit the exceptional and simple case $\alpha=1$).

 It was clear that the integrals of the
 form (0.2), (0.14)  have no expression
  in terms of classical special functions,
 but they were important for probabilists and attracted their
 interest, see \cite{Fel0}, \cite{Fel}, \cite{Zol1}.
 The basic text on this subject is Zolotar\"ev's book \cite{Zol2},
 1986,
 see also bibliography in this book.

 Levy also introduced stable stochastic processes
 (see \cite{Lev3}). Non-explicitness of
  stable densities make stable processes difficult for investigations;
  nevertheless some  collection
  of explicit formulae is known,
  see Dynkin \cite{Dyn},
  Neretin \cite{Ner}, Pitman, Yor \cite{PY}.

In this paper, the expression (0.15)
 appears in the formulae (0.10), (1.10).
Also
the formulae (0.11),
(0.6) are variants of
the "multiplication theorem for the stable
laws "\cite{Zol2}, Theorem 3.3.1; there are many other places,
where we touch formulae from Zolotarev' book \cite{Zol2},
I do not try to fix all similarities in formulae.

\smallskip

5) The functions $\L_{\alpha,\beta}$ arise in a relatively
natural way in the theory of the Laplace transform
(the 'operation calculus'), see below.
The tables of McLachlan, Humbert, Poli \cite{MHP},
1950, contain 18 partial cases of
the integral transformations
defined below; also the transformations (0.10) are contained
in Zolotarev \cite{Zol1},
and a similar construction with Wright functions
is a subject of
 Agarwal \cite{Aga}.

\smallskip

6) It is known (see \cite{MR}) that  pseudodifferential
equations with constant coefficients of the type
$$\Bigl(\sum_{k=0}^n a_k D_{k\alpha}\Bigr) \,f=0$$
admit explicit analysis. Apparently this phenomena
is related
to the identities (0.7), (0.12).

\smallskip

{\bf 0.4. Structure of the paper.}
In Section 1, we discuss various definitions of the functions
$\L_{\alpha,\beta}$, theirs integral representations,
and also some integrals containing products of two
functions $\L_{\alpha,\beta}$.

In Section 2, we discuss the space $\mathcal K$
 of holomorphic functions defined above; also
 we introduce a standard scale $H_\mu$ of
 Hilbert spaces of holomorphic
 functions in a half-plane. The latter spaces
 are well-known in representation theory of $SL_2(\R)$.

 For our purposes, the space $\mathcal K$ and the
 Hardy space $\H^2$ are almost sufficient.

 In Section 3, we introduce a simple
  construction in a spirit of
  the Vilenkin--Klimyk book \cite{VK}.
  We consider the 6-dimensional
 solvable Lie group of operators
 $$
 f(x)\mapsto\lambda x^h f(ax^\alpha);
 \qquad \text{$\lambda\in\C^*$, $h\in\C,a>0$,
 $\alpha\in\R\setminus 0$}
 $$
 acting in a space of functions
  on a half-line and consider the image of
  this group under the Laplace transform.
  As a result, we obtain a group of
  continuous operators, whose kernels
  are expressed in terms of $\L_{\alpha,\beta}$.
  The most interesting property of these operators
  is the identity (0.7) given above.

In Section 4, we consider a similar construction.
We start from a 7-dimensional semigroup
(4,2) of operators
acting in a space of holomorphic function
on half-plane
and consider its image under the
inverse Laplace transform.
As a result, we obtain a semigroup of integral operators
acting in an appropiate
 space of functions on half-line.

\medskip

{\large \bf 1. Some properties of the functions $\L_{\alpha,\beta}$.}

 \addtocounter{sec}{1}
\setcounter{equation}{0}

\medskip

{\bf 1.1. Definition.}
We define the function $\L_{\alpha,\beta}$ as
the Barnes integral
\begin{equation}
\L_{\alpha,\beta}(z)=
\frac 1{2\pi i} \int_{-i\infty}^{+i\infty}
\Gamma(s)\Gamma(\beta-\alpha s)\,z^{-s} ds
\end{equation}
We must explain a meaning of elements of this formula.

1) Our indices are in the domain $\alpha\in\R$, $\beta\in\C$. Assume
also
\begin{equation}
\beta+\alpha m+ n\ne 0\qquad \text{for all $n,m=0,1,2,\dots$}
\end{equation}

2) Our integral is convergent if
$|\arg z|<(1+\alpha)\pi/2$.

3) Our integrand has poles at the points
$$
s=0,-1,-2,\dots \qquad \text{and}
   \qquad s=\beta/\alpha, (\beta+1)/\alpha, (\beta+2)/\alpha,\dots
$$

Now, we consider two cases $\alpha>0$ and $\alpha<0$.

First, let $\alpha>0$. For $\beta>0$,
we can assume that the contour
of the integration is the imaginary axis
$i\R$ and we leave the pole $s=0$
on the left side from the contour
(denote this contour by $+0+i\R$).
 Otherwise, we consider a contour
$L$ coinciding with $i\R$ near $\pm i\infty$
and separating the left series of poles
($s=-n$) and the right series
$s=(\beta+n)/\alpha$ of poles.
Such contour exists due the condition (1.2).

We also can transform this contour integral
to
$$
\int_L=\int_{+0+i\R} -
\sum_{n:\,(\beta+n)/\alpha<0}
\res_{s=(\beta+n)}
$$

Second, let $\alpha<0$.
Then we consider an arbitrary contour $L$
coinciding with the imaginary axis near $\pm\infty$
and leaving all the poles of the integrand on the left side.
If  $\beta>0$, then we can
choice $L$ being $+0+i\R$.

\smallskip

{\sc Remark.}
For fixed $\beta$, $z>0$, the function $\L_{\alpha,\beta}(z)$
as a function of the parameter $\alpha$ is $C^\infty$-smooth at
$\alpha=0$%
 \footnote{See the integral representation (1.9),
 we  differentiate it in $\alpha$ and
 apply the Lebesgue dominant convergence theorem.}
but it is not real analytic in $\alpha$ at this point
(compare (1.4) and (1.5)). Thus it is not quite clear,
is natural to consider $\L_{\alpha,\beta}$
 as one function or as two functions defined for $\alpha>0$
and $\alpha<0$. For local purposes of this paper,  the first
variant is more convenient.

\smallskip

{\bf 1.2. Expansion of $\L_{\alpha,\beta}$ into  power series.}
We write expansions of $\L_{\alpha,\beta}$ in series applying
the standard Barnes method, see \cite{Sla}, \cite{Mar}.

a) Let $0<\alpha<1$. Then the integral (1.1) is
the sum of residues at the   points $s=-n$, i.e.,
\begin{equation}
\L_{\alpha,\beta}(z)=
\sum_{n=0}^\infty \frac{(-1)^n\Gamma(\alpha n+\beta)} {n!}
    z^n
\end{equation}
This function is well defined   on the whole
complex plane $z\in\C$.
Due (1.2), the $\Gamma$-functions in numerators
 have no poles.

b) Let $\alpha>1$. Then (1.1) is the sum of residues at
the points
$s=(\beta+n)/\alpha$, i.e.,
\begin{equation}
\L_{\alpha,\beta}(z)=
-
\frac 1\alpha
\sum_{n=0}^\infty \frac{(-1)^n\Gamma\bigl((n+\beta)/\alpha\bigr)}
 {n!}z^{(n+\beta)/\alpha}
\end{equation}
Here we assume $z^\nu=\exp(\nu\ln z)$ and
$\ln z\in \R$ for $z>0$.
The series is convergent in the domain $|\arg z|<\infty$
(i.e., our function is defined on the
 universal covering
           surface
$\wt\C^*$ of the punctured complex plane $\C^*=\C\setminus 0$).

\smallskip

c) For $\alpha<0$, the integral (1.1) is the sum
of residues at all the poles,
i.e.,
\begin{equation}
\L_{\alpha,\beta}(z)=
\sum_{n=0}^\infty \frac{(-1)^n\Gamma(\alpha n+\beta\bigr)} {n!}
    z^n           -
\frac 1\alpha
\sum_{n=0}^\infty \frac{(-1)^n\Gamma\bigl((n+\beta)/\alpha)}
 {n!}z^{(n+\beta)/\alpha}
\end{equation}
This expression is valid if poles are simple, i.e.,
$\beta+\alpha n+m\ne 0$ for $n,m=0,1,2,3,\dots$.
But the points  $\beta=-\alpha n-m$ are not really
singular, in these cases some of poles of the integrand
have order 2, and we must apply a formula
for a residue in a non-simple pole
(or remove singularities in (1.5)).

The domain of convergence of (1.5) is
$|\arg z|<\infty$.

\smallskip

{\bf 1.3. A symmetry.}

\smallskip

{\sc Lemma 1.1.}
a) {\it For
 $\alpha>0$, $\beta>0$,}
\begin{equation}
\L_{\alpha,\beta}(z)=
\alpha^{-1} z^{-\beta/\alpha}
\L_{1/\alpha,\beta/\alpha}(z^{-1/\alpha})
\end{equation}

b) {\it For $0<\alpha<1$,}
$$
\L_{\alpha,\alpha}(z)=-\frac 1{\alpha z}[\L_{\alpha,1}(z)-1]
$$

c) {\it For $\alpha>0$,}
\begin{equation}
\L_{\alpha,1}(z)=1-\L_{1/\alpha,1}(z^{-1/\alpha})
\end{equation}

{\sc Proof.}
a) Substituting $t=\beta-\alpha s$ to (1.1), we obtain
\begin{equation}
\L_{\alpha,\beta}(z)=
\frac 1 {2\pi i\alpha}\int_{-\i\infty}^{+i\infty}
\Gamma\bigl((\beta-t)/\alpha\bigr)\Gamma(t)\,z^{-(\beta-t)/\alpha}
dt
\end{equation}
as it was required.

Statement b) follows from
$$
\L_{\alpha,\alpha}(z)=
\sum_{n=0}^\infty
\frac{ (-1)^n\Gamma(\alpha n+\alpha)}
  {n!} z^n=
\frac 1\alpha
\sum_{n=0}^\infty
\frac{ (-1)^n\Gamma(\alpha n+\alpha)(\alpha n+\alpha)}
  {(n+1)!} z^n
$$

c) Substitute $\beta=1$ to (1.6) and
assume $\alpha>1$. Applying b),
we obtain the required statement
for $\alpha>1$.
But the identity (1.6) is symmetric with respect
to the transformation $\alpha\mapsto 1/\alpha$,
$z\mapsto z^{-1/\alpha}$.

{\sc Remark.} The statement c)
 is a well-known symmetry in the theory of
stable distributions, see \cite{Fel}, (17.6.10), see also
\cite{Zol2}, Section 2.3.

\smallskip

{\bf 1.4. Some integral representations of $\L_{\alpha,\beta}$.}

{\sc Lemma 1.2.} {\it Let
 $\alpha\in\R$, $\Re u>0$, $\Re v>0$, $\Re h>0$. Then}
\begin{equation}
\int_0^\infty x^{h-1}\exp(-u x^\alpha-v x)\, dx=
v^{-h}\L_{\alpha,h}(u/v^\alpha)
\end{equation}

{\sc Proof.}
It is very easy to verify this for $\alpha>0$.

1) for $0<\alpha<1$:
we expand the factor $\exp(-u x^\alpha)$ in (1.9)
in Taylor series and integrate term-wise
$$
\sum_{n=0}^\infty \frac{(-u)^n}{n!}
\int_0^\infty x^{\alpha n+h-1} e^{-vx}dx
=\sum_{n=0}^\infty \frac{(-u)^n}{n!}
\cdot \frac{\Gamma(\alpha n+h)}
{v^{\alpha n+h}}
$$

2) Similarly, for $\alpha>1$, we expand the factor
$\exp(-vx)$ into a Taylor series
\begin{multline*}
\sum_{n=0}^\infty \frac{(-v)^n}{n!}
\int_0^\infty
x^{h+n-1} \exp(-ux^\alpha)\,dx
=\\=
\frac 1\alpha
\sum_{n=0}^\infty \frac{(-u)^n}{n!}\cdot
\frac{\Gamma((h+n)/\alpha)}{u^{(h+n)/\alpha}}
=\frac 1\alpha u^{-h/\alpha}\L_{1/\alpha,h/\alpha}(vu^{-1/\alpha})
\end{multline*}
and apply the symmetry (1.6)

3) The case $\alpha<0$ is not obvious, and
we give a calculation that is valid for all
$\alpha\in\R$.
Consider the space $L^2$ on $\R_+$ with respect to the measure
$dx/x$. The left hand side of (1.9) is the $L^2$-inner product of
the functions $\Phi_1$, $\Phi_2$  given by
$$
\Phi_1(x)=\exp(-u x^\alpha);\qquad
\Phi_2(x)=x^{\ov h}\exp(-\ov v x)
$$
The Mellin transform of $\Phi_1$ is
$$
\wt \Phi_1 (\lambda)=
\int_0^\infty x^{\lambda-1}\exp(-u x^\alpha)\,dx=
\frac {\sgn(\alpha)}\alpha  \int_0^\infty
\exp(-uy)y^{\lambda/\alpha-1}\,dy=
\frac {\sgn(\alpha)\Gamma(\lambda/\alpha) }
{\alpha u^{\lambda/\alpha}}
$$
The Mellin transform of $\Phi_2$ is
$$
\wt \Phi_2 (\lambda)=
\int_0^\infty x^{\lambda+\ov h-1}\exp(-\ov v x)\,dx =
v^{-\ov h+\lambda}\Gamma(\ov h)
$$
By the Plancherel formula for the Mellin transform,
we have
$$\int_0^\infty \Phi_1(x)\ov{\Phi_2(x)} \,\frac{dx}x
=\frac 1{2\pi}
\int_{-\infty}^\infty \wt\Phi_1(is) \ov{\wt\Phi_2(is)}
\,ds
$$
i.e.,
$$
\int_0^\infty x^{h-1}\exp(-u x^\alpha-v x)\, dx=
\frac{\sgn(\alpha)}{2\pi\alpha v^h}
\int_{-\infty}^{+\infty}
\Gamma(is/\alpha)\Gamma(h+is)\,
u^{-is/\alpha}  v^{is}\,ds
$$
Then we introduce the new variable $t=s/\alpha$.
\hfill $\boxtimes$

\smallskip

{\bf 1.5. Integral representations. Variants.}
Now let $x$, $y>0$. Let $\alpha<1$, $\theta>0$. Then
\begin{multline}
\frac 1{2 i}
\int_{-i\infty}^{+i\infty}
p^{\theta-1}\exp(-p^\alpha x +py)\, dp
=-\Im\bigl[ y^{-\theta} e^{\pi i\theta}
   \L_{\alpha,\theta} (xy^{-\alpha} e^{\pi i\alpha})\bigr]
=\\=
-\sum_{n=0}^\infty \frac{(-1)^n\Gamma(n\alpha+\theta)}{n!}
   x^n y^{-n\alpha-\theta} \sin(n\alpha+\theta)
\end{multline}
where the integration is given over the imaginary axis,
and $p^\alpha$ is positive real for $p>0$.

\smallskip

Indeed, we  represent the integral in the form
\begin{multline*}
 \frac 1{2 i}e^{\pi i \theta/2}\int_0^\infty
   \exp(-t^\alpha xe^{\pi i\alpha/2} +ity)\,t^{\theta-1} dt
-\\-
\frac 1{2 i}
e^{-\pi i\theta/2}\int_0^\infty
   \exp(-t^\alpha xe^{-\pi i\alpha/2} -ity)\,t^{\theta-1} dt
\end{multline*}
and apply (1.9) with $u=x\exp(\pm\pi i\alpha/2)$,
$v=y\exp(\mp\pi i/2)$.

\smallskip

{\sc Remark.} For $0<\alpha<1$ the expression (1.10) is
a  density of a stable subordinator.

\smallskip

{\sc Remark.} The calculation given above survives
 for the case $\alpha<0$,
 The factor $\exp(-t^\alpha x e^{\pi i\alpha/2})$
 is flat at $t=0$ and hence we can consider
  $\theta<0$. We
transform
$$
\Gamma(\alpha n+\theta)\sin[(\alpha n+\theta)\pi]
=
\pi/\Gamma(1-\theta-\alpha n)
$$
and we reduce (1.10) to the form
$
y^{-\theta}\W_{-\alpha,1-\theta}(x/y^\alpha)
$,
where $\W$ is the Wright function.

\smallskip

{\bf 1.6. Remark. An $L^2(\R)$-inner product.}
Let $\alpha>0$,
$\beta>0$. Let $x$, $y>0$. Consider the function
\begin{equation}
\Psi_{\alpha,\beta,y}(x):=x^{\beta-1}\exp(-yx^\alpha)
\end{equation}
By (1.9), its Laplace transform is
\begin{equation}
\wht
\Psi_{\alpha,\beta,y} (\xi)=
\int_0^\infty x^{\beta-1} \exp\bigl\{-yx^\alpha-\xi x\bigr\}\,dx
=
\xi^{-\beta} \L_{\alpha,\beta}(y/\xi^\alpha)
\end{equation}
Evaluating the $L^2(\R)$-inner
product of $\Psi_{\alpha_1,\beta_1,y_1}$ and
$\Psi_{\alpha_2,\beta_2,y_2}$, we obtain
\begin{align}
\int_0^\infty x^{\beta_1-1}\exp(-y_1x^{\alpha_1})
x^{\beta_2-1}\exp(-y_2 x^{\alpha_2})\,dx\nonumber=\\
=\alpha_2^{-1}\,y_2^{-(\beta_1+\beta_2-1)/\alpha_2}
\L_{\alpha_1/\alpha_2, (\beta_1+\beta_2-1)/\alpha_2}
(y_1 y_2^{-\alpha_1/\alpha_2})
\end{align}
(we substitute $t=x^{\alpha_2}$ and apply Lemma 1.2).
The expression in the right-hand side is symmetric
with respect to
$(\alpha_1,\beta_1)\leftrightarrow(\alpha_2,\beta_2)$
by  Lemma 1.1.

By the Plancherel formula for the Fourier transform,
the same expression can be written in the form
\begin{equation}
\frac 1{2\pi}\int_{-\infty}^{+\infty}
(it)^{-\beta_1} \L_{\alpha_1,\beta_1}(y_1/(it)^{\alpha_1})
(-it)^{-\beta_2} \L_{\alpha_2,\beta_2}(y_2/(-it)^{\alpha_2})
\,dt
\end{equation}

Thus, the expression (1.14) equals  (1.13).

\smallskip

{\sc Remark.}
The integral
(1.14) looks like a kernel of a product of two integral operators;
moreover (1.13) shows that this product has the same form, i.e.,
we obtain a family of integral operators
closed with respect to multiplication.
Below we propose two ways to
give a precise sense for this observation;
apparently, there are other possibilities.

\smallskip

{\bf 1.7. Remark. Convolutions.}
Preserve notation (1.11), (1.12).
We have
$$\Psi_{\alpha,\beta,y}(x) \Psi_{\alpha,\beta',y'}(x)=
\Psi_{\alpha,\beta+\beta',y+y'}(x)
$$
Hence
$$
\frac 1{2\pi i}
\int_{-i\infty}^{i\infty}
\wht \Psi_{\alpha,\beta,y}(u)
\wht\Psi_{\alpha',\beta',y'}(z-u)du
=\wht\Psi_{\alpha,\beta+\beta',y+y'}(z)
$$

\medskip

{\bf\large 2. Spaces of holomorphic functions. Preliminaries}

\nopagebreak

 \addtocounter{sec}{1}
\setcounter{equation}{0}

\medskip

{\bf 2.1. Spaces $L^2_\mu(\R_+)$.}
Fix $\mu>0$. Denote by $L^2_\mu(\R_+)$
the space $L^2$ on the half-line $\R_+$, $x>0$,
with respect to the weight
$\Gamma(\mu)^{-1}x^{\mu-1}dx$, i.e., the Hilbert space
with the inner product
$$
\langle f,g\rangle_{[\mu]}:=
\frac 1{\Gamma(\mu)}
\int_0^\infty f(x)\ov{g(x)}x^{\mu-1} dx
$$

For instance,
\begin{equation}
\langle \exp(-z x),\exp(-u x)\rangle_{[\mu]}=
(z+\ov u)^{-\mu}
\end{equation}
for arbitrary complex $u$, $z$ satisfying
$\Re z>0$, $\Re u>0$.

\smallskip

{\bf 2.2. Hilbert spaces of holomorphic functions
on a half-plane.}
Let $\Pi$ be the right half-plane $\Re z>0$
on the complex plane.
We consider the Hardy space $\H^2$ on
$\Pi$. Recall that this space consists of functions
holomorphic in the half-plane, whose boundary values
on the
imaginary axis
$\Re z=0$ exist and are contained in $L^2(\R)$.

The inner product in $\H^2$ is given by
\begin{align*}
\langle f,g\rangle=
\frac 1{2\pi}\int_{-\infty}^{+\infty} f(it)\ov{g(it)}\,dt
:=
\lim_{\epsilon\to+0}
\frac 1{2\pi}\int_{-\infty}^{+\infty}
f(\epsilon+it)\ov{g(\epsilon+it)}\,dt
\end{align*}

The Hardy space is an element of the following
 one-parametric scale
$H_\mu$, $\mu>0$, of spaces of holomorphic functions.

Fix $\mu>1$.
Consider the space $H_\mu=H_\mu(\Pi)$
consisting of functions $f(z)$ holomorphic
in $\Pi$ and satisfying the condition
$$
\int_\Pi |f(z)|^2 (\Re z)^{\mu-2} dz \, d\ov z<\infty
$$
where $dz \, d\ov z$ denotes the Lebesgue measure on $\Pi$.

We define an inner product  in $H_\mu$
by
\begin{equation}
\langle f,g \rangle_\mu=
\frac{\mu-1}\pi
\int_\Pi f(z) \ov{g(z)} (\Re z)^{\mu-2} dz \, d\ov z
\end{equation}
The space $H_\mu$ is a Hilbert space with respect
 to this inner  product.

The reproducing kernel%
\footnote{For machinery of reproducing kernels, see,
for instance, \cite{Ber}, \cite{Ner2}.}
 of this space is
$$
K_\mu(z,u)=(z+\ov u)^{-\mu}
$$
This means that the function $\Xi_u(z)$
given by
$$\Xi_u(z)=(z+\ov u)^{-\mu}$$
satisfies the reproducing property
\begin{equation}
\langle f,\Xi_u\rangle_\mu= f(u)\qquad \text{for all $f\in H_\mu$.}
\end{equation}
In particular,
\begin{equation}
\langle \Xi_u, \Xi_w\rangle_\mu=(w+\ov u)^{-\mu}
\end{equation}

The  space $H_\mu$ can be defined by (2.3), (2.4)
without reference to explicit formula (2.2) for the inner product.
Indeed, consider an abstract Hilbert space $H$ with a system of vectors
$\Xi_u$, where $u\in\Pi$, and assume that their inner
products age given by
(2.4). Such Hilbert space exists, see formula (2.1).
Assume also that linear combinations of $\Xi_u$ are dense in $H$.
Then for each $h\in H$ we consider the holomorphic function
on $\Pi$ given by
\begin{equation}
f_h(u):=\langle h,\Xi_u\rangle_H
\end{equation}
and thus we have identified our space $H$
 with some space of holomorphic functions on $\Pi$.

 \smallskip

But the last construction survives for arbitrary $\mu>0$
(since the existence of $H_\mu$ is provided by formula (2.1)
and this formula is valid for $\mu>0$).

\smallskip

{\sc Remark.} For $\mu=1$ we obtain the Hardy space
$\H^2$.

\smallskip

{\sc Remark.}
For $0<\mu<1$, it is possible to write an integral
formula for the inner product in $H_\mu$ involving
derivatives. But it is more convenient to use
the definition (2.3)--(2.4)
 or to consider the analytic continuation
of the integral (2.2) with respect to $\mu$.

\smallskip

We define the {\it weighted Laplace transform} $\LL_\mu$
by
$$
\LL_\mu f(z)
=\frac 1{\Gamma(\mu)}
\int_0^\infty f(x)\exp(-ux)\,x^{\mu-1} dx
$$
For $\mu=1$ we obtain the usual Laplace transform
$\LL=\LL_1$.
The following statement is well-known%
\footnote{The case $\mu=1$ is a Paley--Wiener theorem.}.

\smallskip

{\sc Lemma 2.1.}
{\it The weighted Laplace transform
is a unitary operator}
$$\LL_\mu:\,\,\,L^2_\mu(\R_+)\to H_\mu(\Pi)$$

{\sc Proof.} Consider a function
$$\xi_u(x)=\exp(-ux)$$
in $L^2_\mu$. Its image under $\LL_\mu$
is $\Xi_u$. It remains to compare (2.1) and (2.4).

 \smallskip

{\sc Remark.} The transform $\LL_\mu$ is precisely the
operator defined by the formula (2.5). Indeed, we can assume
$H=L_\mu$, then the corresponding space of functions $f_h$
is $H_\mu$.

\smallskip

{\bf 2.3. Operators in the spaces $H_\mu$.}
Recall a standard general trick,
apparently discovered by Berezin \cite{Ber} (formulae
(2.6), (2.9) given below are  valid in
arbitrary Hilbert space defined by
a reproducing kernel).

Let $A$ be a {\it bounded} operator $H_\mu\to H_\mu$.
Define the function
\begin{equation}
M(z,u)=A \Xi_u(z)
\end{equation}
Then $M(z,u)$ is the {\it kernel} of the operator $A$.
For $\mu>1$ we can write literally
\begin{equation}
Af(z)=
\frac{\mu-1}\pi
\int_\Pi M(z,u)f(u) (\Re u)^{\mu-2}\,du\,d\ov u
\end{equation}
Respectively, for $\mu=1$,
\begin{equation}
Af(z)=\frac 1{2\pi}
\int_{-\infty}^\infty
M(z,-it) f(it)\,dt
\end{equation}
For general $\mu>0$, we can write
\begin{equation}
A f(z)=\langle f, \ov{M(z,u)}\rangle_\mu
\end{equation}
the the inner product is given in the space
$H_\mu$ of functions depending in the variable $u$.
Formulae (2.7)--(2.8) are partial cases of this formula.
In particular, the integrals (2.7)--(2.8) are convergent for
each $z\in\Pi$ and $f\in H_\mu$.

{\bf 2.4. Space of rapidly decreasing functions.}
We also consider the  space $\mathcal K=\H^2\cap \mathcal S(i\R)$, where
$\mathcal S(i\R)$
is the Schwartz space (consisting of functions
on the imaginary axis rapidly decreasing
with all derivatives).

We  can say that $\mathcal K$ is the space of functions holomorphic in
$\Re z>0$ and continuous in $\Re z\ge 0$ such that

a) $f(it)$ is $C^\infty$-smooth

b) For each $k>0$ and $N>0$ there exists $M$ such that
\begin{equation}
|f^{(k)}(z)|\le M (1+|z|)^{-N}
\end{equation}

Consider the space
$\mathcal S(\R_+)$ consisting of smooth functions $f$ on
$[0,\infty)$ such that

a) $f^{(k)}(0)=0$ for all $k\ge 0$

b) $\lim_{x\to+\infty} f^{(k)}(x) x^N=0$ for all $N>0$, $k\ge 0$.

In other words, $\mathcal S(\R_+)$ is the intersection
of the Schwartz space $\mathcal S(\R)$ on $\R$ and the space
$L^2(\R_+)$.

\smallskip

{\sc Lemma 2.2.} {\it
The space $\mathcal K$ is the image of $\mathcal S(\R_+)$
 under the Laplace transform.}

 \smallskip

{\sc Proof.} Let $f\in \mathcal S(\R_+)$. Integrating by parts, we obtain
$$
\int_0^\infty f^{(k)}(x)e^{-px}dx=
p^k\int_0^\infty f(x)e^{-px}dx
$$
The left-hand side is a bounded function in $p$,
looking to the right-hand side, we observe that $(\LL f) (p)$
is rapidly decreases for $\Re p\ge 0$.

 Conversely, a function $F$ satisfying (2.10)
 is an element of $\H^2$. Hence $f=\LL^{-1}F$ is supported
 by $\R_+$. Since (2.10) is valid for
 $z\in i\R$, we have $f\in \mathcal S(\R)$.

{\bf 2.5. Fractional derivations.}
We define the operators of
{\it fractional differentiation}
 $D_h$ in 
  $\mathcal K$
by
\begin{equation}
D_h f(z)=
\frac{\Gamma(h+1)}{2\pi }
\int_{-\infty}^{+\infty}
\frac{f(it)\,dt}{( -it+z)^{h+1}}
\end{equation}
A branch of $\theta(z,it)=( -it+z)^{h+1}$
is determined from the condition $\theta(x,0)>0$ for $x>0$.

\smallskip

{\sc Lemma 2.3.}
a) {\it $D_h$ is an operator $\mathcal K \to \mathcal K$
for each $h\in \C$. }

\smallskip

b) {\it For integer $n>0$,}
$$
D_n f(z)= (-1)^n\frac {d^n}{dz^n} f(z)
$$

c) {\it For positive integer $m$,}
\begin{align}
D_{-m} f(z):=
\lim_{s\to m}
\frac{\Gamma(-s+1)}{2\pi i}
\int_{-\infty}^{+\infty}
(-it+z)^{s-1}f(t)\,du
=\\
=(-1)^m
\int_{-i\infty}^{z} dz_1
\int_{-i\infty}^{z_1}dz_2\dots
\int_{-i\infty}^{z_{n-1}}
f(z_n)\, dz_n\nonumber
\end{align}

d) $D_{h_1} D_{h_2}=D_{h_1+h_2}$

{\sc Proof.} a) Convergence of the integral for $\Re u>0$ is obvious.
Let us show rapid decreasing  of $g:=D_h f$ at $z\to \infty$.

Let $\Re h<0$.
We represent  (2.11) as a contour integral
$$
\frac{\Gamma(h+1)}{2\pi i}
\int_{-i\infty}^{+i\infty}
\frac{f(u)\,du}{(-u+z)^{h+1}}
$$
Denote $R:=|z|$.
Then we replace a part $(-iR,+iR)$ of the contour
of the integration $i\R$ by the semi-circle
$R \exp(i\phi)$, where $\phi\in(-\pi,\pi)$.
Since (2.10), all the 3 summands of the integral
rapidly tend
to 0 as $R$ tend to $\infty$.

If $\Re h\ge 0$, we integrate our expression  by parts
and obtain
$$
\frac{(-1)^k\Gamma(h-k+1)}{2\pi i}
\int_{-i\infty}^{+i\infty}
(z-u)^{k-h-1} f^{(k)}(u)\,du
$$
We choose $k>\Re(h+1)$ and repeat the same consideration.

Also,
$$\frac
d{dz} D_h f= -D_{h+1}f
$$
and this implies rapid decreasing (2.10) of derivatives.

\smallskip

b) This is the Cauchy integral representation for derivatives.

\smallskip

c) First, we give a remark that formally is not necessary.
Factor $\Gamma(1+h)$ has a pole for $h=-m$. Let us show that $\int$
vanishes at this point. Indeed, we have the expression
$$
\int_{-i\infty}^{i\infty} (z-u)^{m-1}f(u)\,du
$$
We replace a part $(-iR,+iR)$ of the contour
of the integration $i\R$ by the semi-circle
$R \exp(i\phi)$ as above
and tend $R$ to $\infty$.

Now give a  proof of c). Consider the operator
of  indefinite integration
\begin{equation}
If(z)=\int_{-i\infty}^z f(u)\,du
\end{equation}
Changing the contour as above, we obtain
$If\in\mathcal K$.

For $f\in \mathcal K$ we have
\begin{align*}
\Gamma(1-s)\int_{-i\infty}^{+i\infty} (z-u)^{s-1} f(u)\,du=
\\
=(-1)^m\Gamma(1-s+m)\int_{-i\infty}^{+i\infty} (z-u)^{s-m-1} (I^m f)(u)\,du=
\end{align*}
Now we can substitute
$s=m$.

d) This is valid for $\Re h_1=\Re h_2=0$ by the
statement b) of  following Lemma 2.4.
Then we consider the analytic continuation in $h$.

\smallskip

{\sc Lemma 2.4.} a) {\it For $s\in\R$, the operator $D_{is}$ is a
unitary
operator in each $H_\mu$. Its kernel (in the sense of $H_\mu$)
is                                       }
$$
\frac {\Gamma(\mu+is)}{\Gamma(is)} (z+\ov u)^{-\mu-is}
$$

b) $D_{is_1} D_{is_2}=D_{is_1+is_2}$

\smallskip

{\sc Proof.} a) Consider the operator
$U_{i s}$ in $L^2_\mu(\R_+)$ given by
$$U_{i s} f(x)=f(x) \,x^{i s}$$
Let us evaluate the kernel of the operator
$$
\LL_\mu U_{is} \LL^{-1}_\mu: H_\mu\to H_\mu
$$
By (2.6), we must evaluate the function
$\LL_\mu U_{is} \LL^{-1}_\mu  \Xi_u$.
We have
\begin{align*}
\LL^{-1}_\mu  \Xi_u(x)&=\exp(-\ov u x)
,\\
U_{is}\LL^{-1}_\mu  \Xi_u(x)&=\exp(-\ov u x) x^{is}
;\\
\LL_\mu U_{is} \LL^{-1}_\mu  \Xi_u(z)&=(z+\ov u)^{-\mu-is}
{\Gamma(\mu+is)}/{\Gamma(is)}
\end{align*}

For $\mu=1$ our operator coincides
with the operator $D_{is}$ defined above.
In fact, all the operators
$\LL_\mu U_{is} \LL^{-1}_\mu: H_\mu\to H_\mu$
induce the same operator
$D_{is}$ in  $\mathcal K$.
Indeed, the operator
$L^{-1}_\mu L_\nu$ is the operator
of multiplication
by $x^{\mu-\nu}$, and this operator commutes with
$U_{is}$.

\smallskip

b)  corresponds to the identity
$x^{i(s_1+s_2)}=x^{is_1} x^{is_2}$
after the Laplace transform.

\medskip

{\bf\large 3. Operators in  spaces of holomorphic functions}

 \addtocounter{sec}{1}
\setcounter{equation}{0}

\medskip

{\bf 3.1. Some operators acting in $\mathcal S(\R_+)$.}
We consider the following one-parameter groups of operators
in the space $\mathcal S(\R_+)$
(see  2.4).
\begin{align}
U_\alpha f(x)=f(x^\alpha),\qquad \alpha>0;
\\
V_a(x)  f(x)=f(ax),\qquad a>0;
\\
W_h f(x)=x^h f(x),\qquad h\in \C
\end{align}
The last group is a complex one-parameter
group, i.e., a real two-parameter group.
The infinitesimal generators of these groups are
respectively
$$
E_1 f(x)= x\ln x \frac d{dx} f(x);
\qquad
E_2 f(x)=x \frac d{dx} f(x);
\qquad E_3 f(x)= (\ln x ) f(x)
$$
They satisfy the commutation relations
\begin{equation}
[E_1,E_2]=-E_2; \qquad [E_1,E_3]=E_3; \qquad
[E_2,E_3]=1
\end{equation}
Thus we obtain a real 6-dimensional Lie algebra
$\frak g$
spanned by the operators
$$
E_1,\quad E_2,\quad  E_3,\quad  iE_3,\quad  1,\quad i
$$
The algebra  $\frak g$ is solvable and it contains a two-dimensional
center $\R\cdot 1+ \R\cdot i$.

Also $\frak g$ is a real subalgebra (but not a
real form)  in 4-dimensional
complex Lie algebra
$$\C\cdot E_1+ \C\cdot E_2 +\C\cdot E_3+\C\cdot 1$$

Obviously,
\begin{align}
U_\alpha W_h U_\alpha^{-1}= W_{\alpha h}
          ;\\
U_\alpha V_a U_\alpha^{-1}= V_{a^{1/\alpha}}
\end{align}

Consider the group $G$ generated by the one-parameter
groups (3.1)--(3.3).
General element of this group is an operator of the form
$$\lambda\cdot R(h,\alpha,a)$$
where $\lambda\in\C^*$ and
\begin{equation}
R(h,\alpha,a) f(x)= x^{h} f(a x^\alpha)
=W_h U_\alpha V_a f(x)
\end{equation}
The product is given by
\begin{equation}
R(g,\beta,b) R(h,\alpha,a)=
b^h R(g+h\beta,\alpha \beta,ab^\alpha)
\end{equation}

We also can add the operator
$$U_{-1} f(x)=f(1/x)$$
or equivalently we can allow $\alpha\in\R\setminus 0$
in (3.1), (3.7), (3.8).  Then we obtain a Lie group
$G^\circ$ of operators
containing two connected components; the group $G$ defined above
is the connected component containing 1.

\smallskip

{\bf 3.2. Operators in $L^2_\mu(\R_+)$.}
Now, let us fix $\mu>0$.
Then the operators
$$|\alpha|^{1/2} a^{1/2} R((\alpha-\mu)/2+is,\alpha,a)$$
 are unitary in
$L^2_\mu(\R_+)$.
Such operators form a 4-dimensional solvable Lie group
with an one-dimensional center; denote this group
by $G_\mu^\circ$.

\smallskip

{\sc Remark.} All other operators $R(h,\alpha,a)$
are unbounded in $L^2_\mu(\R)$.

\smallskip

{\bf 3.3. Operators in $H_\mu$.}
Let us evaluate the kernel in $H_\mu$ of the operator
$\LL_\mu R(h,\alpha,a) \LL_\mu^{-1}$
using 2.3.

We have
\begin{align}
\LL_\mu^{-1} \Xi_u (x)&=\exp(-\ov u x); \nonumber\\
R(h,\alpha,a) L_\mu^{-1} \Xi_u(x)&=
x^h \exp(- a \ov u x^\alpha);\nonumber\\
\LL_\mu R(h,\alpha,a) \LL_\mu^{-1} \Xi_u(z)&=
\frac 1{\Gamma(\mu)}
\int_0^\infty  x^{h+\mu-1} \exp(- a \ov u x^\alpha- zx)\,dx
=\\&=
z^{-\mu-h}\L_{\alpha,h+\mu}(\ov u a/z^\alpha)
\end{align}
and we obtain (for $\mu>1$) the integral operator
\begin{equation}
\wt R(h,\alpha,a)F(z) =
\frac{1}{\pi\Gamma(\mu-1)}
z^{-\mu-h}
\int_{\Pi}
\L_{\alpha,h+\mu}(\ov u a/z^\alpha)
F(u)\,(\Re u)^{\mu-2} du\,d\ov u
\end{equation}

For $\mu=1$ we understand   (3.11) as (2.8),
and for $\mu<1$ as (2.9).

Formally, our algorithm of evaluating of the kernel
is valid only for bounded operators.
Thus we proved the following theorem

\smallskip

{\sc Theorem 3.1.} {\it Let
$$\Re h=(\alpha-\mu)/2$$
Then an operator $\wt R(h,\alpha,a)$
defined by (3.11) is
 unitary in $H_\mu$ up to a scalar factor.
The product of two
operators $\wt R$ is  given by the formula
\begin{equation}
\wt R(g,\beta,b)\wt R(h,\alpha,a)=
b^h \wt R(g+h\beta,\alpha \beta,ab^\alpha)
\end{equation}
These operators generate a
4-dimensional solvable Lie group isomorphic
to the group $G_\mu^\circ$ described in 3.2.}

\smallskip

{\bf 3.4. Operators in the space $\mathcal K$.}
The group $G^\circ$ defined in 3.1 acts in
the space $\mathcal S(\R_+)$.
Since the weighted Laplace transform
 $\LL_\mu$ identifies $\mathcal S(\R_+)$
and $\mathcal K$, this group also acts in $\mathcal K$.
For the subgroup $G_\mu^\circ\subset G^\circ$, the action
was constructed in the previous subsection, but
 formula (3.11) for the kernel
almost survives for a  general element of $G^\circ$.
We consider separately $\alpha>0$ and $\alpha<0$.

\smallskip

a) Let $\alpha>0$.
We substitute $y=x^\alpha$ to the expression (3.9)
and obtain
\begin{equation}
\frac 1 {|\alpha|}
\int_0^\infty
y^{(h+\mu)/\alpha-1} \exp(-a\ov u y-z y^{1/\alpha})\,dy
\end{equation}

If $h$ satisfies the condition
\begin{equation}
(\Re h+\mu)/\alpha>0
\end{equation}
then our integral is convergent
(otherwise, we have a non-integrable
singularity at 0).
 The expression (3.13) is the Laplace transform of
the function
$y^{(h+\mu)/\alpha-1} \exp(-z y^{1/\alpha})$.
Since this function is an element of $L^1$,
its Laplace transform is a bounded function in $\Pi$.
Hence, for $F\in\mathcal K$, the integral (3.11) is convergent
and is holomorphic in $h$. Thus the formula (3.11)
defines the $\LL_\mu$-image of $R(h,\alpha,a)$
for all triples $h$, $\alpha$, $a$ satisfying (3.14).

Next, we write
\begin{equation}
R(h,\alpha,a)=W_{-n}\circ R(h+n,\alpha,a)
\end{equation}
For sufficiently large $n$, the $\LL_\mu$-image
of $R(h+n,\alpha,a)$
is defined by formula (3.11); the  $\LL_\mu$-image
of $W_{-n}$ is the iterated indefinite integration (2.12).

b) $\alpha<0$. We again transform the integral to the form
(3.13). Now the integrand is smooth at 0; for convergence at infinity
we are need in the condition
$(\Re h+\mu)/\alpha<0$.
Then we repeat (3.15)%
\footnote{This is not really necessary, since for $\Re u>0$
integral (3.13) is convergent. But for the case $\mu=1$
it is pleasant to have an expression for $\Re u=0$. }

    \smallskip

Thus we obtain the following theorem

\smallskip

{\sc Theorem 3.2.}
{\it Fix $\mu>0$.
 Let $h\in\C$, $\alpha\in\R\setminus 0$, and $a>0$.
 If $\Re h+\mu>0$, then we define the integral operator
 $\wt R(h,\alpha,a)$
in $\mathcal K$ by (3.11). Otherwise,
we consider $n$ such that $\Re h+n>0$ and define
$$
\wt R(h,\alpha,a)=(-1)^n I^n\circ \wt R(h+n,\alpha,a)
$$
there $I$ is the operator of indefinite
integration (2.13)
Then all these operators
 are bounded in $\mathcal K$ and their product
is given by (3.12). The group generated by $\wt R(h,\alpha,a)$
is isomorphic to the group $G^\circ$ defined in 3.1.}

\smallskip

{\sc Remark.} For different $\mu$ we obtain
 the same group of operators  in $\mathcal K$;
but identification of this group with $G^\circ$
depends on $\mu$.

\smallskip

{\bf 3.5. Statements formulated in Introduction.}
Let $\mu=1$. The operator $A_\alpha$ given by (0.5)
is the $\LL_1$-image  of $U_\alpha$, the fractional
derivation $D_h$ is the $\LL_1$-image  of $W_h$,
and the operator $R_a$ is the  $\LL_1$-image  of
$V_a$. Now (0.7), (0.9) follow from (3.5), (3.6).

The operator $F\mapsto zF$ is the $\LL$-image of $d/dx$,
and this implies (0.8).

\smallskip

{\bf 3.6. Hankel type transforms.}
The kernel of the operator $\wt R(h,-1,1)$ is
$$
K(z,u)=\mathrm{const}\cdot \int_0^\infty
 x^{h+\mu-1}\exp (-\ov u/x-z x)\,dx
$$
This expression is a modified Bessel function of Macdonald, see
\cite{HTF}. The corresponding
integral transform is similar to the Hankel transform.

\smallskip

{\bf 3.7. Another group of symmetries.}
Now we consider the group of unitary operators in $L^2(\R_+)$
generated by
\begin{align*}
U'_{\alpha}f(x)&=
|\alpha|^{1/2} f(x^\alpha) x^{(\alpha-1)/2};\\
V'_a f(x)&=a^{1/2}f(ax);\\
T_\beta(is) f(x)&=\exp(is x^\beta) f(x),
 \qquad s\in \R, \beta in\R
\end{align*}
This group is infinite dimensional since
it contains all the operators having the form
$$f(x)\mapsto f(x)\exp(i \sum_{j=1}^N s_j x^{\beta_j} )$$
Obviously, we have
\begin{equation}
U_\alpha T_\beta(is)  U_\alpha^{-1} f(x)=T_{\alpha\beta}(is);
\qquad
V_a T_\beta(is)   V_a^{-1} =T_\beta(is a^\beta)
\end{equation}

Consider the image of our group of operators under the standard
Laplace transform $\LL$. The operators $U'_\alpha$, $V'_a$
are contained in the group $G^\circ$ and their images
$\wt U'_\alpha$, $\wt V'_a$  are described above.
The image of $T_\beta(is)$ is the convolution operator
$$
\wt  T_\beta(is) F(z)=\frac 1{2\pi i}
\int_{-i\infty}^{i\infty}
M(z-u) F(u)du$$
where
$$
M(z)=\int_0^\infty \exp(is x^\beta -zx)\,dx=
\L_{\beta,1}(is/z^\beta)
$$
In particular, we obtain the identity
$$
\wt U_\alpha \wt T_\beta(is) \wt  U_\alpha^{-1} f(x)=
\wt T_{\alpha\beta}(is)
$$
for our operators with $\L$-kernels.

\medskip

{\bf\large 4. Operators in space of functions on half-line}

\medskip

\nopagebreak

 \addtocounter{sec}{1}
\setcounter{equation}{0}

{\bf 4.1. Spaces of functions.}
Consider the space $\P$ consisting of $C^\infty$-functions
on half-line $x\ge 0$ satisfying

\smallskip

a) $f^{(k)}(0)=0$ for all $k\ge 0$

\smallskip

b) $\lim_{x\to+\infty} f^({m})(x) \exp(-\epsilon x)=0$ for all $\epsilon>0$, $m\ge 0$.

\smallskip

Consider also
the space $\F$, whose elements are functions holomorphic in
the half-plane $\Re z>0$ satisfying the condition

---  for each $\epsilon>0$, $N>0$ where exists $C$
such that
$$
 |f(z)|< C /|z|^{-N} \qquad \text{for $\Re z>\epsilon$}
$$

{\sc Lemma 4.1.}
{\it The image of the space $\P$ under the Laplace transform
is $\mathcal F$.}

{\sc Proof.} a) For $f\in\mathcal P$ denote
$F=\LL f$. Fix small $\delta>0$.
$$
F(p)=\int_0^\infty f(x)e^{-px}dx=
\int_0^\infty \bigl[f(x)e^{-\delta x}\bigr] e^{-(p-\delta)x}dx
$$
Since $f(x)e^{-\delta x}$ is an element of the Schwartz space,
for each $N$ we have an estimate
$$
|F(p)|\le C|p-\delta|^{-N}
$$
for $\Re p\ge\delta$. For $\Re p>2\delta$, we can write
$$
|F(p)|\le  2^N C |p|^{-N}
$$
Thus $F\in \mathcal F$.

b) Let $F\in \mathcal F$.
The inversion formula for $\LL$ gives
$$
f(x)=\frac 1{2\pi i}\int_{a-i\infty}^{a+i\infty}
e^{px}F(p)\,dp=
\frac 1{2\pi} e^{xa}\int_{-i\infty}^{i\infty} e^{itx} F(a+it)\,dt
$$
Since $F(a+it)$ is an element of Schwartz space $\mathcal S(\R)$,
the function
$$e^{-ax}f(x)=
\frac 1{2\pi} \int_{-i\infty}^{i\infty} e^{itx} F(a+it)\,dt
$$
is an element of $\mathcal S(\R)$. By the Paley--Wiener theorem
this function is supported by $\R^+$.
Since $a>0$ is arbitrary, we obtain we required statement.

\smallskip

{\bf 4.2. Fractional derivations.}
Consider the Riemann--Liouville operators
 of fractional integration (see \cite{SKM})
\begin{equation}
J_r f(x)=\frac 1{\Gamma(r)}
\int_0^x {f(y)}(x-y)^{r-1}\, dy
\end{equation}
in the space $\P$. The integral is convergent for $r>0$.
For integer positive $r=n$ we have
$$
J_r f(x)=\int_0^x dx_1\int_0^{x_1}dx_2\dots
 \int_0^{x_{n-1}} f(x_n)\,dx_n
 $$
For fixed $f$ and $x$, the function
$J_r f(x)$ admits a holomorphic continuation to
the whole plane $r\in \C$, and for integer
negative $r=-n$ we have
$$
J_{-n} f(x) =\frac {d^n} {dx^n} f(x)
$$
We also have the identity
$$
J_r J_p=J_{r+p}
$$
for all $r$, $p\in\C$.

Laplace transform $\LL$ identifies the Riemann--Liouville
fractional integrations $J_r$
with the operators in $\F$ given by
$$F(z)\mapsto z^{-r} F(z)$$

{\bf 4.3. Operators.}  Now we consider the semigroup $\Gamma$
of operators in $\F$ consisting of  transformations
$$
\lambda Q(\theta,\alpha,a)
$$
where $\lambda \in \C^*$,
\begin{equation}
 Q(\theta,\alpha,a) F(z) =z^\theta F(a z^\alpha)
\end{equation}
and the parameters satisfy the conditions
\begin{equation}
0<\alpha<1,\qquad
\arg a+\alpha\pi/2<\pi/2,\qquad \arg a-\alpha \pi/2>-\pi/2,
\qquad \theta\in \C
\end{equation}

{\sc Remark.}
The restrictions (4.4) mean that $z\in \Pi$
implies $az^\alpha\in\Pi$.

Obviously, we have
\begin{equation}
Q(\theta',\alpha',a') Q(\theta,\alpha,a)
=(a')^\theta Q(\theta'+\theta\alpha',\alpha\alpha',a (a')^\alpha)
\end{equation}
Thus $\Gamma$ is a 7-dimensional semigroup with 2-dimensional
center.

\smallskip

{\sc Remark.} The semigroup $\Gamma$  can be  embedded to
a 7-dimensional Lie group.
The parameters of this group are $\lambda\in\C^*$, $a\in\C^*$,
$\theta\in \C$, $a>0$.
 The multiplication in this group
is determined by the formula (4.4),
where $\alpha>0$ and $\theta$, $a\in \C$.
But  the corresponding operators (4.2) are not well-defined
in the space $\F$.

\smallskip

The $\LL^{-1}$-image of the operators $Q(\theta,\alpha,a)$
is given by the formula
\begin{equation}
\wt  Q(\theta,\alpha,a) f(x)=
\int_0^\infty N(x,y) f(y)\, dy
\end{equation}
where the kernels $N(x,y)$ are given by
$$
N(x,y)=\frac 1{2\pi}
\int_{-i\infty}^{i\infty}
z^\theta\exp(-a z^\alpha y+ zx)\,dz
$$
For $\theta>-1$,
 we transform this integral in the same way as in 1.5
and obtain
\begin{multline*}
N(x,y)=\frac 1 {2\pi i y^{\theta+1}}
\Bigl\{
e^{i(\theta+1)\pi/2} \L_{\alpha,\theta+1}
(y a x^{-\alpha} e^{i\pi \alpha/2})
-\\-
e^{-i(\theta+1)\pi/2} \L_{\alpha,\theta+1}
(y a x^{-\alpha} e^{-i\pi \alpha/2})
\Bigr\}
\end{multline*}
For real $a>0$, $\theta>-1$
 we have an expression of the form (1.10).

If $\theta<-1$, we write
$z^\theta F(az^\alpha)=z^{-n}z^{\theta+n} F(az^\alpha)$
and
\begin{equation}
\wt Q(\theta,\alpha,a)
:=J_n\circ \wt Q(\theta+n,\alpha,a)
\end{equation}
for sufficiently large $n$.

\smallskip

{\sc Theorem 4.1.} {\it For $\alpha$, $a$
satisfying (4.3), the operators $\wt Q(\theta,\alpha,a)$
defined by (4.5), (4.6) are bounded in
$\mathcal P$ and their product is given by (4.4).
}

\smallskip

Formulae (0.11)--(0.12) given in Introduction
are particular cases of this statement.

{\sf Math.Phys. Group,
Institute of Theoretical and Experimental Physics,

B.Cheremushkinskaya, 25, Moscow 117259

\& University of Vienna, Math. Dept.,
Nordbergstrasse, 15, Vienna 1090, Austria

neretin@mccme.ru

\end{document}